\documentclass{gtart} 

\input gtoutput
\volumenumber{3}\papernumber{16}\volumeyear{1999}
\pagenumbers{397}{404}\published{30 November 1999}
\proposed{Joan Birman}\seconded{Shigeyuki Morita, Dieter Kotschick}
\received{21 July 1999}
\accepted{23 November 1999}

\newtheorem{thm}{Theorem}[section]   
\newtheorem{lem}[thm]{Lemma}         
\theoremstyle{definition}
\newtheorem{defn}[thm]{Definition}   
\newtheorem*{rem}{Remark}            
\newcommand{\Z}{{\bf{Z}}}

\begin{document}

\title{The Burau representation is not faithful for $n = 5$}
\shorttitle{The Burau representation is not faithful for n = 5}
\asciititle{The Burau representation is not faithful for n = 5}
\authors{Stephen Bigelow}                  
\address{Department of Mathematics\\UC Berkeley\\Berkeley, CA 94720, USA} 
\email{bigelow@math.berkeley.edu}

\begin{abstract} 
The Burau representation is a natural
action of the braid group $B_n$
on the free $\Z[t,t^{-1}]$--module of rank $n-1$.
It is a longstanding open problem
to determine for which values of $n$ this representation is faithful.
It is known to be faithful for $n=3$.
Moody has shown that it is not faithful for $n\geq 9$
and Long and Paton improved on Moody's techniques
to bring this down to $n\geq 6$.
Their construction uses a simple closed curve
on the $6$--punctured disc
with certain homological properties.
In this paper we give such a curve on the $5$--punctured disc,
thus proving that the Burau representation is not faithful for $n\geq 5$.
\end{abstract}

\asciiabstract{The Burau representation is a natural
action of the braid group B_n
on the free Z[t,t^{-1}]-module of rank n-1.
It is a longstanding open problem
to determine for which values of n this representation is faithful.
It is known to be faithful for n=3.
Moody has shown that it is not faithful for n>8
and Long and Paton improved on Moody's techniques
to bring this down to n>5.
Their construction uses a simple closed curve
on the 6-punctured disc
with certain homological properties.
In this paper we give such a curve on the 5-punctured disc,
thus proving that the Burau representation is not faithful for n>4.}

\primaryclass{20F36}                
\secondaryclass{57M07, 20C99}              
\keywords{Braid group, Burau representation}

\maketitlepage

\section{Introduction}
\label{INTRODUCTION}

The braid groups $B_n$ appear in many different guises.
We recall here the definition we will be using
and some of the main properties we will need.
For other equivalent definitions
see \cite{birman:blam}.

Let $D$ denote a disc
and let $q_1,\dots,q_n$ be $n$ distinct points
in the interior of $D$.
For concreteness, 
take $D$ to be the disc in the complex plane 
centered at the origin and having radius $n+1$,
and take $q_1,\dots,q_n$ to be the points $1,\dots,n$.
Let $D_n$ denote the punctured disc $D\setminus \{q_1,\dots,q_n\}$,
with basepoint $p_0$ on $\partial D$, say $p_0 = -(n+1)i$.

\begin{defn}
The braid group $B_n$
is the group of all equivalence classes of 
orientation preserving homeomorphisms $h\co D_n\rightarrow D_n$
which fix $\partial D$ pointwise,
where two such homeomorphisms are equivalent 
if they are homotopic rel $\partial D$.
\end{defn}

It can be shown that
$B_n$ is generated by $\sigma_1,\dots,\sigma_{n-1}$,
where $\sigma_i$ exchanges punctures $q_i$ and $q_{i+1}$ 
by means of a clockwise twist.

Let $x_1,\dots,x_n$ be free generators of $\pi_1(D_n,p_0)$,
where $x_i$ passes counterclockwise around $q_i$.
Consider the map $\epsilon\co \pi_1(D_n) \rightarrow \Z$
which takes a word in $x_1,\dots,x_n$
to the sum of its exponents.
Let $\tilde{D}_n$ be the corresponding covering space.
The group of covering transformations of $\tilde{D}_n$ is $\Z$,
which we write as a multiplicative group generated by $t$.
Let $\Lambda$ denote the ring $\Z[t,t^{-1}]$.
The homology group $H_1(\tilde{D}_n)$
can be considered as a $\Lambda$--module,
in which case it becomes a free module of rank $n-1$.

Let $\psi$ be an autohomeomorphism of $D_n$
representing an element of $B_n$.
This can be lifted to a map
$\tilde{\psi}\co \tilde{D}_n\rightarrow\tilde{D}_n$
which fixes the fiber over $p_0$ pointwise.
This in turn induces a $\Lambda$--module automorphism
$\tilde{\psi}_*$ of $H_1(\tilde{D}_n)$.
The {\em (reduced) Burau representation} is the map
\[
\psi \mapsto \tilde{\psi}_*.
\]
This is an $(n-1)$--dimensional representation of $B_n$ over $\Lambda$.

The main result of this paper is the following.

\begin{thm}
\label{burau5}
The Burau representation is not faithful for $n = 5$.
\end{thm}

The idea is to use the fact that
the Dehn twists about two simple closed curves
commute if and only if those simple closed curves can be
freely homotoped off each other.
Our construction will use two simple closed curves 
which cannot be freely homotoped off each other
but in some sense ``fool'' the Burau representation 
into thinking that they can.
To make this precise, we first make the following definition.

\begin{defn}
Suppose $\alpha$ and $\beta$ are two arcs in $D_n$.
Let $\tilde{\alpha}$ and $\tilde{\beta}$ be lifts
of $\alpha$ and $\beta$ respectively
to $\tilde{D}_n$.
We define
\[
\int_\beta \alpha= \sum_{k\in\Z} (t^k\tilde{\alpha},\tilde{\beta}) t^k,
\]
where $(t^k\tilde{\alpha},\tilde{\beta})$
denotes the algebraic intersection number
of the two arcs in $\tilde{D}_n$.
Note that this is only defined 
up to multiplication by a power of $t$,
depending on the choice of lifts $\tilde{\alpha}$ and $\tilde{\beta}$.
This will not pose a problem because we will only be interested in
whether or not $\int_\beta \alpha$ is zero.
\end{defn}

\begin{thm}
\label{lp}
For $n\geq 3$, the Burau representation of $B_n$ is not faithful
if and only if there exist embedded arcs $\alpha$ and $\beta$ on $D_n$
such that $\alpha$ goes from $q_1$ to $q_2$,
$\beta$ goes from $p_0$ to $q_3$ or from $q_3$ to $q_4$,
$\alpha$ cannot be homotoped off $\beta$ rel endpoints,
and $\int_\beta \alpha = 0$.
\end{thm}

The special case in which $\beta$ goes from $p_0$ to $q_3$
follows easily from \cite[Theorem 1.5]{long-paton:burau}.
This special case is all we will need to prove Theorem~\ref{burau5}.
Nevertheless, 
we will give a direct proof of Theorem~\ref{lp}
in Section~\ref{LP}.
In Section~\ref{BURAU5} we give a pair of curves
on the $5$--punctured disc
which satisfy the requirements of Theorem~\ref{lp},
thus proving Theorem~\ref{burau5}.

Throughout this paper,
elements of the braid group act on the left.
If $\psi_1$ and $\psi_2$ are elements of the braid group $B_n$
then we denote their commutator by:
\[
[\psi_1,\psi_2] = \psi_1^{-1} \psi_2^{-1} \psi_1 \psi_2.
\]

\section{Proof of Theorem \ref{lp}}
\label{LP}

It will be useful to keep the following lemma in mind.
It can be found in \cite[Proposition 3.10]{flp:tdt}.

\begin{lem}
\label{no_digons}
Suppose $\alpha$ and $\beta$ are simple closed curves on a surface
which intersect transversely at finitely many points.
Then $\alpha$ and $\beta$ can be freely homotoped
to simple closed curves which intersect at fewer points
if and only if there exists a ``digon'',
that is, an embedded disc
whose boundary consists of one subarc of $\alpha$
and one subarc of $\beta$.
\end{lem}

First we prove the ``only if'' direction of Theorem~\ref{lp}.
Let $n\geq 3$ be such that
for any embedded arcs $\alpha$ from $q_1$ to $q_2$
and $\beta$ from $p_0$ to $q_3$ in $D_n$
satisfying $\int_\beta \alpha = 0$
we have that $\alpha$ can be homotoped off $\beta$ rel endpoints.
Let $\psi\co D_n \rightarrow D_n$ 
lie in the kernel of the Burau representation.
We must show that $\psi$ is homotopic to the identity map.

Let $\alpha$ be the straight arc from $q_1$ to $q_2$
and let $\beta$ be the straight arc from $p_0$ to $q_3$.
Then $\int_\beta \psi(\alpha) = 0$.
Thus $\psi(\alpha)$ can be homotoped off $\beta$.
By applying this same argument to an appropriate conjugate of $\psi$
we see that $\psi(\alpha)$ can be homotoped off
the straight arc from $p_0$ to $q_j$ for any $j=3,\dots,n$.
It follows that we can homotope $\psi$ so as to fix $\alpha$.
Similarly, we can homotope $\psi$ so as to fix
every straight arc from $q_j$ to $q_{j+1}$ for $j=1,\dots,n-1$.
The only braids with this property are powers of
$\Delta$, the Dehn twist about a simple closed curve
parallel to $\partial D$.
But $\Delta$ acts as multiplication by $t^n$ on $H_1(\tilde{D}_n)$.
Thus the only power of $\Delta$ which lies in the kernel of
the Burau representation is the identity.

We now prove the converse for the case in which
$\beta$ is an embedded arc from $q_3$ to $q_4$ in $D_n$.
Let $\alpha$ be an embedded arc from $q_1$ to $q_2$
such that $\alpha$ cannot be homotoped off $\beta$ rel endpoints
but $\int_\beta \alpha = 0$.
Let $\tau_\alpha\co D_n\rightarrow D_n$ be a ``half Dehn twist''
about the boundary of a regular neighborhood of $\alpha$.
This is the homeomorphism which exchanges punctures $q_1$ and $q_2$
and whose square is a full Dehn twist
about the boundary of a regular neighborhood of $\alpha$.
Similarly, let $\tau_\beta$ be a half Dehn twist
about the boundary of a regular neighborhood of $\beta$.
We will show that the commutator of $\tau_\alpha$ and $\tau_\beta$
is a non-trivial element of the kernel of the Burau representation.

\begin{figure}[ht!]\small
  \centering

\begin{picture}(90,25)

\thinlines

\put(10,10){\oval(20,20)[l]}
\put(10,20){\vector(1,0){10}}
\put(20,22){$\alpha'$}
\put(20,20){\line(1,0){10}}
\qbezier(30,20)(40,20)(50,10)
\qbezier(50,10)(60,0)(70,0)
\put(70,0){\line(1,0){20}}
\put(90,10){\oval(20,20)[r]}
\put(90,20){\vector(-1,0){10}}
\put(80,20){\line(-1,0){10}}
\qbezier(70,20)(60,20)(50,10)
\qbezier(50,10)(40,0)(30,0)
\put(30,0){\line(-1,0){20}}

\thicklines
\put(10,10){\line(1,0){80}}
\put(10,10){\circle*{4}}
\put(90,10){\circle*{4}}

\end{picture}

  \caption{The ``figure eight'' $\alpha'$}
  \label{fig_eight}
\end{figure}
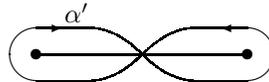

Let $\epsilon$ be an embedded arc in $D_n$ which crosses $\alpha$ once.
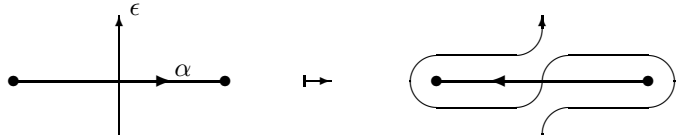
\begin{figure}[ht!]\small
  \centering

\begin{picture}(270,45)

\thinlines

\put(40,0){\vector(0,1){45}}
\put(43,45){$\epsilon$}

\thicklines

\put(0,20){\line(1,0){80}}
\put(60,20){\vector(1,0){0}}
\put(60,22){$\alpha$}
\put(0,20){\circle*{4}}
\put(80,20){\circle*{4}}

\thinlines

\put(110,20){\vector(1,0){10}}
\put(110,18){\line(0,1){4}}

\put(210,0){\oval(20,20)[tl]}
\put(210,10){\line(1,0){30}}
\put(240,20){\oval(20,20)[r]}
\put(240,30){\line(-1,0){30}}
\put(210,20){\oval(20,20)[tl]}
\put(190,20){\oval(20,20)[br]}
\put(190,10){\line(-1,0){30}}
\put(160,20){\oval(20,20)[l]}
\put(160,30){\line(1,0){30}}
\put(190,40){\oval(20,20)[br]}
\put(200,40){\vector(0,1){5}}

\thicklines

\put(160,20){\line(1,0){80}}
\put(180,20){\vector(-1,0){0}}
\put(160,20){\circle*{4}}
\put(240,20){\circle*{4}}

\end{picture}

  \caption{The action of $\tau_\alpha$}
  \label{fig_twist}
\end{figure}
Figure~\ref{fig_twist} shows $\epsilon$ and
its image under the action of $\tau_\alpha$.
Thus the effect of $\tau_\alpha$
on $\epsilon$ is, up to homotopy rel endpoints,
to insert the ``figure-eight'' $\alpha'$
shown in Figure~\ref{fig_eight}.
Now let $\tilde{\epsilon}$ be a lift of $\epsilon$ 
to the covering space $\tilde{D}_n$.
Note that $\alpha'$ lifts to a closed curve in $\tilde{D}_n$.
Thus the effect of $\tilde{\tau}_\alpha$ 
on $\tilde{\epsilon}$ is, up to homotopy rel endpoints,
to insert a lift of~$\alpha'$.

Let $\tilde{\gamma}$ be a closed arc in $\tilde{D}_n$.
The effect of $\tilde{\tau}_\alpha$ on $\tilde{\gamma}$
is to insert some lifts of $\alpha'$.
If we consider $\tilde{\gamma}$ and $\tilde{\alpha}'$
as representing elements of $H_1(\tilde{D}_n)$ then
\[
(\tilde{\tau}_\alpha)_*(\tilde{\gamma}) 
= \tilde{\gamma} + P(t)\tilde{\alpha}',
\]
where $P(t) \in \Lambda$.
Similarly,
\[
(\tilde{\tau}_\beta)_*(\tilde{\gamma}) 
= \tilde{\gamma} + Q(t)\tilde{\beta}',
\]
where $Q(t) \in \Lambda$
and $\beta'$ is a figure eight
defined similarly to $\alpha'$.

Any lift of $\alpha$, and hence of $\alpha'$,
has algebraic intersection number zero
with any lift of $\beta$.
It follows that 
\[
(\tilde{\tau}_\beta)_*(\tilde{\alpha}') = \tilde{\alpha}'.
\]
Thus
\[
(\tilde{\tau}_\beta \tilde{\tau}_\alpha)_*(\tilde{\gamma})
= (\tilde{\gamma} + Q(t)\tilde{\beta}') + P(t)\tilde{\alpha}'.
\]
Similarly
\[ 
(\tilde{\tau}_\alpha \tilde{\tau}_\beta)_*(\tilde{\gamma})
= (\tilde{\gamma} + P(t)\tilde{\alpha}') + Q(t)\tilde{\beta}'.
\]
Thus
$(\tilde{\tau}_\alpha)_*$ and $(\tilde{\tau}_\beta)_*$ commute, 
so the commutator $[\tau_\alpha, \tau_\beta]$ 
lies in the kernel of the Burau representation.

It remains to show that $[\tau_\alpha, \tau_\beta]$ 
is not homotopic to the identity map.
Let $\gamma$ be the boundary of a regular neighborhood of $\alpha$.
Using Lemma~\ref{no_digons} and the fact that
$\gamma$ cannot be freely homotoped off $\beta$,
it is not hard to check that $\tau_\beta(\gamma)$ 
cannot be freely homotoped off $\alpha$.
A similar check then shows that $\tau_\alpha \tau_\beta(\gamma)$
cannot be freely homotoped off $\tau_\beta(\gamma)$.
Thus $\tau_\alpha \tau_\beta(\gamma)$ is not freely homotopic to
$\tau_\beta(\gamma)$.
But $\tau_\beta(\gamma) = \tau_\beta \tau_\alpha(\gamma)$.
Thus $\tau_\alpha \tau_\beta$ is not homotopic to
$\tau_\beta \tau_\alpha$,
so $[\tau_\alpha, \tau_\beta]$ is not homotopic to the identity map.

The case in which $\beta$ goes from $p_0$ to $q_3$
can be proved by a similar argument.
Instead of a half Dehn twist about 
the boundary of a regular neighborhood of $\beta$
we use a full Dehn twist 
about the boundary of a regular neighborhood of 
$\beta \cup \partial D$.
Instead of a figure eight curve $\beta'$
we obtain a slightly more complicated curve
which is a commutator of $\partial D$
and the boundary of a regular neighborhood of $\beta$.

\section{Proof of Theorem \ref{burau5}}
\label{BURAU5}

\begin{figure}[ht!]\small
  \centering

\unitlength0.7pt

\begin{picture}(340,240)(4,0)

\thinlines

\put(324,0){\vector(0,1){66}}
\put(326,66){$\beta$}
\put(324,66){\line(0,1){66}}

\put(252,132){\oval(24,48)[t]}
\put(246,132){\oval(60,72)[t]}
\put(246,132){\oval(84,96)[t]}
\put(246,132){\oval(108,120)[t]}
\put(252,132){\oval(144,144)[t]}

\put(228,132){\oval(24,48)[b]}
\put(234,132){\oval(60,72)[b]}
\put(234,132){\oval(84,96)[b]}
\put(234,132){\oval(108,120)[b]}
\put(234,132){\oval(132,144)[b]}

\put(96,132){\oval(24,48)[t]}
\put(102,132){\oval(60,72)[t]}
\put(102,132){\oval(84,96)[t]}
\put(102,132){\oval(108,120)[t]}
\put(102,132){\oval(132,144)[t]}

\put(96,132){\oval(120,144)[b]}
\put(54,132){\oval(12,36)[b]}
\put(108,132){\oval(48,48)[b]}
\put(108,132){\oval(72,72)[b]}
\put(108,120){\line(0,1){12}}

\thicklines

\put(18,144){\line(1,0){78}}
\put(18,136){\line(1,0){312}}
\put(18,128){\line(1,0){312}}
\put(228,120){\line(1,0){78}}

\put(18,150){\oval(12,12)[bl]}
\put(18,150){\oval(28,28)[bl]}
\put(18,114){\oval(28,28)[tl]}
\put(330,142){\oval(12,12)[br]}
\put(330,142){\oval(28,28)[br]}
\put(306,114){\oval(12,12)[tr]}

\put(4,150){\line(0,1){72}}
\put(12,150){\line(0,1){72}}
\put(18,222){\oval(12,12)[tl]}
\put(18,222){\oval(28,28)[tl]}

\put(18,228){\line(1,0){312}}
\put(18,236){\line(1,0){312}}

\put(330,222){\oval(12,12)[tr]}
\put(330,222){\oval(28,28)[tr]}

\put(336,142){\line(0,1){80}}
\put(344,142){\line(0,1){80}}

\put(4,42){\line(0,1){72}}
\put(18,42){\oval(28,28)[bl]}
\put(18,28){\line(1,0){140}}
\put(298,28){\vector(-1,0){140}}
\put(155,33){$\alpha$}
\put(298,42){\oval(28,28)[br]}
\put(312,42){\line(0,1){72}}

\put(54,120){\circle*{4}}
\put(96,144){\circle*{4}}
\put(108,120){\circle*{4}}
\put(228,120){\circle*{4}}
\put(252,144){\circle*{4}}

\end{picture}

  \caption{Arcs on the $5$--punctured disc}
  \label{fig_b5}
\end{figure}
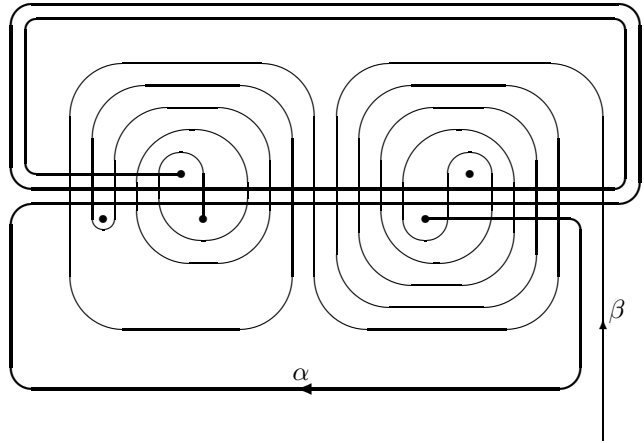

Let $\alpha$ and $\beta$ be the embedded arcs on $D_5$
as shown in Figure~\ref{fig_b5}.
These cannot be homotoped off each other rel endpoints,
as can be seen by applying Lemma \ref{no_digons}
to boundaries of regular neighborhoods of $\alpha$ and $\beta\cup\partial D$.
It remains to show that $\int_{\beta}\alpha = 0$.

Let $\tilde{\alpha}$ and $\tilde{\beta}$
be lifts of $\alpha$ and $\beta$ to $\tilde{D}_5$.
Each point $p$ at which $\beta$ crosses $\alpha$
contributes a monomial $\pm t^k$ to $\int_\beta \alpha$.
The exponent $k$ is such that 
$\tilde{\beta}$ and $t^k \tilde{\alpha}$ cross at a lift of $p$,
and the sign of the monomial is the sign of that crossing.
We choose our lifts and sign conventions such that
the first point at which $\beta$ crosses $\alpha$
is assigned the monomial $+t^0$.

In Figure \ref{fig_b5},
the sign of the monomial at a crossing $p$ will be
positive if $\beta$ is directed upwards at $p$
and negative if $\beta$ is directed downwards at $p$.
The exponents of the monomials can be computed 
using the following remark:

\begin{rem}
Let $p_1,p_2 \in \alpha\cap\beta$
and let $k_1$ and $k_2$ be the exponents of the monomials
at $p_1$ and $p_2$ respectively.
Let $\alpha'$ and $\beta'$ be the arcs from $p_1$ to $p_2$ 
along $\alpha$ and $\beta$ respectively
and suppose that $\alpha' \cap \beta' = \{p_1,p_2\}$.
Let $k$ be such that
$\alpha' \cup \beta'$ bounds a $k$--punctured disc.
Then $|k_2 - k_1| = k$.
If $\beta'$ is directed counterclockwise around the $k$--punctured disc
then $k_2 \geq k_1$,
otherwise $k_2 \le k_1$.
\end{rem}

One can now progress along $\beta$,
using the above remark to calculate the exponent at each crossing
from the exponent at the previous crossing.
Reading the crossings from left to right, top to bottom,
we obtain the following:
\begin{eqnarray*}
\int_{\alpha} \beta 
&=& - t^{-3} - t^0 + t^1 + t^{-1} + t^{-3} \\
& &\mbox{}-t^{-1}-t^2+t^3 + t^1 +t^{-1}-t^{-2}- t^0 - t^2 + t^1 +t^{-2} \\
& &\mbox{}-t^{-1}+t^0-t^1 + t^2 - t^3 + t^2 - t^1 + t^0 - t^{-1} + t^{-2} \\
& &\mbox{}- t^1 - t^4+t^5 + t^3 + t^1 - t^0 - t^2 - t^4 + t^3 + t^0 \\
& &\mbox{}- t^1 + t^2-t^3 + t^4 - t^5 + t^4 - t^3 + t^2 - t^1 + t^0 \\
& &\mbox{}- t^2 + t^1-t^0 + t^{-1} - t^{-2} \\
&=& 0. \\
\end{eqnarray*}
Thus $\alpha$ and $\beta$ satisfy the requirements of Theorem~\ref{lp},
and we conclude that the Burau representation is not faithful for $n=5$.

The proof of Theorem~\ref{lp}
gives an explicit non-trivial element of the kernel,
namely the commutator of a half Dehn twist
about the boundary of a regular neighborhood of $\alpha$
and a full Dehn twist about the boundary of a regular neighborhood
of $\beta\cup\partial D$.
The following element of $B_5$ sends $\alpha$ to 
a straight arc from $q_4$ to $q_5$:
\[
\psi_1 = 
\sigma_3^{-1} \sigma_2 \sigma_1^2 \sigma_2 \sigma_4^3 \sigma_3 \sigma_2.
\]
The following element of $B_5$ 
sends $\beta$ to a straight arc from $p_0$ to $q_5$:
\[
\psi_2 =
\sigma_4^{-1} \sigma_3 \sigma_2 \sigma_1^{-2} \sigma_2
\sigma_1^2 \sigma_2^2 \sigma_1 \sigma_4^5.
\]
Thus the required kernel element is:
\[
[
\psi_1^{-1}
  \sigma_4 
\psi_1,
\psi_2^{-1} 
  \sigma_4 \sigma_3 \sigma_2 \sigma_1^2 \sigma_2 \sigma_3 \sigma_4 
\psi_2
].
\]
This is a word of length $120$ in the generators.

The arcs in Figure~\ref{fig_b5}
were found using a computer search,
although they are simple enough to check by hand.
A similar computer search for the case $n=4$
has shown that any pair of arcs on $D_4$
satisfying the requirements of Theorem~\ref{lp}
must intersect each other at least $500$ times.

We conclude with an example of a non-trivial braid in 
the kernel of the Burau representation for $n=6$
which is as simple such a braid as one could reasonably hope 
to obtain from Theorem~\ref{lp}.
\begin{figure}[ht!]\small
  \centering

\unitlength0.7pt

\begin{picture}(270,100)

\thinlines

\put(30,45){\oval(60,60)[t]}
\put(135,45){\oval(90,90)[t]}
\put(135,90){\vector(1,0){0}}
\put(135,95){$\beta$}
\put(240,45){\oval(60,60)[t]}

\put(45,45){\oval(90,90)[b]}
\put(60,45){\line(0,-1){15}}
\put(225,45){\oval(90,90)[b]}
\put(210,45){\line(0,-1){15}}

\thicklines

\put(30,45){\vector(1,0){105}}
\put(134,48){$\alpha$}
\put(135,45){\line(1,0){105}}

\put(30,45){\circle*{4}}
\put(60,30){\circle*{4}}
\put(120,60){\circle*{4}}
\put(150,60){\circle*{4}}
\put(210,30){\circle*{4}}
\put(240,45){\circle*{4}}

\end{picture}

  \caption{Arcs on the $6$--punctured disc}
  \label{fig_b6}
\end{figure}
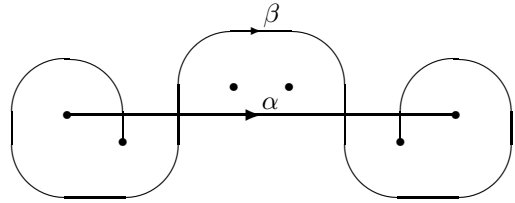
The curves in Figure~\ref{fig_b6}
give us the braid
\[
[\psi_1^{-1}\sigma_3\psi_1, \psi_2^{-1}\sigma_3\psi_2],
\]
where
\[
\psi_1 = \sigma_4 \sigma_5^{-1} \sigma_2^{-1} \sigma_1,
\]
and
\[
\psi_2 = \sigma_4^{-1} \sigma_5^2 \sigma_2 \sigma_1^{-2}.
\]
This is a word of length $44$ in the generators.

\end{document}